\documentclass[12pt]{article}
\usepackage{amsfonts}
\usepackage{graphicx}
\usepackage{a4}
\newcommand\R{{\mathbb{R}}}
\newtheorem{Lemma}{Lemma}
\newtheorem{lemma}[Lemma]{Lemma}

\newtheorem{Theorem}{Theorem}
\newtheorem{Corollary}{Corollary}
\newtheorem{Definition}{Definition}

\newcommand\qed{\hfill\mbox{$\Box$}}

\newcommand\text[1]{\mathrm{#1}}
\author{N. Br\"annstr\"om \ and \ V. Gelfreich\thanks{The authors thank
Prof.~D.~Turaev for useful discussions and helpful suggestions.}
\\[6pt]
{\small
Mathematics Institute,
University of Warwick}\\
{\small
Coventry, CV4 7AL, United Kingdom}\\[6pt]
{\small
\begin{tabular}{ll}
E-mail:&\tt N.L.A.Brannstrom@warwick.ac.uk\\
&\tt  V.Gelfreich@warwick.ac.uk
\end{tabular}
}
}
\title{Drift of slow variables in slow-fast Hamiltonian systems}
\date{October 15, 2007}
\begin{document}
\maketitle

\begin{abstract}
We study the drift of slow variables in a slow-fast Hamiltonian system with
several fast and slow degrees of freedom. For any periodic trajectory
of the fast subsystem with the frozen slow variables we define an action.
For a family of periodic orbits, the action is a scalar function of the slow variables
and can be considered as a Hamiltonian function which generates some slow dynamics.
These dynamics depend on the family of periodic orbits.

Assuming the fast system with the frozen slow variables has a pair of hyperbolic
periodic orbits connected by two transversal heteroclinic trajectories,
we prove that for any path composed of a finite sequence of
slow trajectories generated by action Hamiltonians, there is 
a trajectory of the full system whose slow component
shadows the path.

\end{abstract}

\section{Introduction}

We consider a slow-fast Hamiltonian system described by a smooth Hamiltonian
function 
\(
H(p,q,v,u;\varepsilon)\,.
\)
This system is slow-fast due to a small parameter
in the symplectic form
\[
\Omega=dp\wedge dq+\frac1\varepsilon dv\wedge du\,.
\]
Therefore the equations of motion take the form
\begin{equation}\label{Eq:motion}
\begin{array}{ll}
\displaystyle
\dot q=\frac{\partial H}{\partial p}\,, \qquad&
\displaystyle
\dot p=-\frac{\partial H}{\partial q}\,,\\[10pt]
\displaystyle
\dot u=\varepsilon\frac{\partial H}{\partial v}\,, &
\displaystyle
\dot v=-\varepsilon\frac{\partial H}{\partial u}  \,.
\end{array}
\end{equation}
Equations of this form often arise after rescaling a part of the variables
in a Hamiltonian system with the standard symplectic form.%

The variable $(p,q)$ are called fast and $(v,u)$ are slow.
We assume that the system has $m+d$ degrees of freedom, where $m$
is the number of fast degrees of freedom and $d$ is the number of
slow ones.

After substituting $\varepsilon=0$ into equation (\ref{Eq:motion})
we see that the values of $(v,u)$ remain constant in time
and the system can be interpreted as a family of Hamiltonian
systems with $m$ degrees of freedom which depends on $2d$ parameters.
We call it a {\em frozen system}:
\begin{equation}\label{Eq:frozen}
\begin{array}{ll}
\displaystyle
\dot q=\frac{\partial H}{\partial p}\,, \qquad&
\displaystyle
\dot p=-\frac{\partial H}{\partial q}\,,\\[10pt]
\displaystyle
\dot u=0\,, &
\displaystyle
\dot v=0  \,.
\end{array}
\end{equation}

The case when the fast system has one degree of freedom is relatively well understood.
Indeed, in this case the frozen system typically represents a fast oscillator.
The averaging method can be used to eliminate the dependence on the fast
oscillations from the slow system. Therefore trajectories of the slow system
are close to trajectories of an autonomous system with $2d$ degrees of freedom
over very long time intervals (see e.g. \cite{Arnold78,BM1961}). 

Points of equilibria of the frozen system form surfaces called slow manifolds. 
Normally hyperbolic slow manifolds persists and normally elliptic slow
manifolds do not in general. In both cases the dynamics in a neighbourhood 
of a slow manifold can be described using normal forms (for a discussion
see e.g.~\cite{GelfreichLerman}).

The case when the fast system has more than one degree of freedom is 
notably more difficult. The effect of the fast system on the slow variables
strongly depends on the dynamics of the fast system.
If the frozen fast system oscillates with a constant
vector of frequencies, generalisations of the averaging method can be used
\cite{Lochak,Neishtadt1976}. The averaging method can be also used if the frozen
system is uniformly hyperbolic \cite{Anosov} or, more generally, 
if the frozen system is ergodic
and time averages converge sufficiently fast to space averages~\cite{Kiefer}.
In all these cases the dynamics of the slow variables is described, in the leading
order, by the vector field obtained by taking an average of the slow component
of (\ref{Eq:motion}) over the space of fast variables 
\begin{equation}\label{Eq:average}
\dot u=\varepsilon\left<\frac{\partial H}{\partial v}\right>\,,
\qquad
\dot v=-\varepsilon\left<\frac{\partial H}{\partial u}\right>  \,.
\end{equation}
This approximation strongly relies on the fact that in an ergodic
system the time average over a trajectory equals
the space average for almost all trajectories. 
The approximation error strongly depends on the rate of
convergence for time averages.
If the number of fast degrees of freedom is larger than one,
there is no reason to expect that the time average 
over a periodic orbit converges to the average
over the space. Therefore we should expect that the slow component
of a trajectory whose fast component stays near a periodic orbit
of the frozen system may strongly deviate 
from the averaged dynamics described by (\ref{Eq:average}).
Moreover, we note that periodic orbits are dense in the
case of an Anosov system. 

\medskip

In this paper we assume that the frozen system has a
compact invariant set bearing chaotic dynamics of horseshoe type
created by transversal heteroclinics between two saddle periodic
orbits. This situation typically arises when a periodic orbit has a transversal
homoclinic.

In this invariant set hyperbolic periodic orbits are dense and every two
periodic orbits are connected by a heteroclinic orbit.
We select a finite subset of periodic orbits with relatively short periods.
We construct trajectories of the full system which 
switch between neighbourhoods of the periodic orbits in
a prescribed way. We show that the slow component of such trajectories
drifts in a way quite similar to trajectories of 
a random Hamiltonian dynamical system with $2d$ degrees of freedom. 

The trajectories constructed in this paper strongly deviate from
the averaged dynamics. We think this mechanism is 
responsible for the largest possible rates of deviation. 

A similar construction is used in~\cite{GT2007} for studying
drift of the enrgy in a Hamiltonian system which depends
on time explicitly and slowly. In particular,  
it was shown in~\cite{GT2007} that switching between 
fast periodic orbits indeed provides the fastest rate 
of energy growth in several situations.

\medskip

The rest of the paper has the following structure. 
In Section~\ref{Se:mainthm} we state our main theorem 
and discuss its application to systems with
one slow degree of freedom. In Section~\ref{Se:action}
we describe slow dynamics of the full system~(\ref{Eq:motion}) 
near a family of periodic orbits of the frozen system. 
The description is based on an action associated with the frozen periodic 
orbits and can be of independent interest.
The central ingredient of the proof of the main theorem
is preservation of normally hyperbolic manifolds 
formed by families of uniformly hyperbolic
orbits of the frozen system which is explained in Section~\ref{Se:normalhyperb}.
In this section we explain how symbolic dynamics can be used
to describe the dynamics of the full system restricted to an invariant subset
close to the hyperbolic invariant set of the frozen system.
The discussion is based on ideas of \cite{GT2007}.
Section~5 analyses the long time behaviour of the slow component
of the full dynamics.
The last section of the paper finishes the proof of the main theorem.

\section{Accessibility and drift of slow variables\label{Se:mainthm}}

The total energy is preserved, so we study the dynamics on
a single energy level. Without any loss in generality we may 
consider the dynamics in the zero energy level
\[
{\cal M}_\varepsilon=\{\, H(p,q,v,u;\varepsilon)=0\, \}\,.
\]
First we state our assumptions on the dynamics
of the frozen system. 
Let $D\subset\R^{2d}$ be in a bounded domain. We assume 
\begin{itemize}
\item[{[A1]}] 
the frozen system  has two smooth families
of hyperbolic periodic orbits $L_c(v,u)\subset {\cal M}_0$ defined for all 
$(v,u)\in D$, $c\in\{\,a,b\,\}$. 

\item[{[A2]}] 
the frozen system has two smooth families of transversal heteroclinic orbits:
\begin{eqnarray*}
\Gamma_{ab}(v,u)&\subset& W^u(L_a(v,u))\cap W^s(L_b(v,u))\,,\\
\Gamma_{ba}(v,u)&\subset& W^u(L_b(v,u))\cap W^s(L_a(v,u))\,,\qquad \forall(v,u)\in D\,.
\end{eqnarray*}
\end{itemize}
We note that under these assumptions the frozen system has
a family of uniformly hyperbolic invariant transitive sets $\Lambda_{(v,u)}$,
also known as Smale horseshoes. For every $(v,u)\in D$, this set contains a countable
number of saddle periodic orbits, which are dense in~$\Lambda_{(v,u)}$.
Moreover, every two periodic orbits in~$\Lambda_{(v,u)}$ are
connected by a transversal heteroclinic orbit, which also
belongs to~$\Lambda_{(v,u)}$. It is well known that the dynamics
on the Smale horseshoe can be described using the language
of Symbolic Dynamics. We define
\[
\Lambda:=\bigcup_{(v,u)\in D}\Lambda_{(v,u)}\,.
\]

Before stating our main theorem we give a couple of definitions.
\begin{Definition}
The action of a periodic orbit $L_c$ is defined by the integral
\[
J_c(v,u):=
\oint_{L_c(v,u)}p\,dq\,.
\]
\end{Definition}
The function $J_c$ is independent of the fast variables and can be
considered itself as a Hamiltonian function which generates some slow dynamics
of $(v,u)$ variables:
\begin{equation}\label{Eq:slow}
v'=-\frac1{T_c}\frac{\partial J_c}{\partial u}\,,
\qquad
u'=\frac1{T_c}\frac{\partial J_c}{\partial v}\,,
\end{equation}
where ${}'$ stands for the derivative with respect to the slow
time $\tau=\varepsilon t$, and $T_c$ is the period of $L_c$.
System~(\ref{Eq:slow}) is Hamiltonian 
with the symplectic form 
 $\omega_c=T_c(v,u)dv\wedge du$. Alternatively the equations
can be interpreted as a result of a time scaling in a
standard Hamiltonian system.

In the next sections we will show that for properly chosen initial conditions
the slow component of the corresponding trajectory of (\ref{Eq:motion})
oscillates near a trajectory of this slow Hamiltonian system.

Inside the Smale horseshoe there are infinitely many periodic orbits
connected by transversal heteroclinics. Each periodic orbit
has an action associated with it. 
We select a finite subset of periodic orbits and consider
the collection of their actions.
In general we should expect all those actions to be different. 

In this paper we prove that there are trajectories of the full system
such that their slow components follow any finite path composed of 
segments of slow trajectories generated by actions. 
Those trajectories of the full system shadow a chain composed of the periodic orbits
and heteroclinic trajectories and spend most of the time near
periodic orbits of the frozen system.

Let us give a definition of an accessible path and then state the
theorem. Consider a finite family of functions $J_k:D\to\R$, $k=1,\ldots,n$.
Let $\Phi^\tau_k$ be the Hamiltonian flow with Hamiltonian function $J_k$
and the symplectic form $\omega_k=T_k(v,u)dv\wedge du$ where 
$T_k>0$ is the period of the corresponding orbit.
For every point $z=(v,u)\in D$ we define
\[
\sigma_k(z)=\sup\{\, \tau: \Phi_k^{\tau'}(z)\in D\quad \mbox{for all $\tau'\in (0,\tau)$}\,\}\,,
\]
which is the time required to leave the domain $D$.
If the trajectory is defined for all $\tau>0$ we set $\sigma_k(z)=+\infty$.
Obviously, $\sigma_k(z)>0$ for any $z$ and $k$ due to openness of $D$.

\begin{Definition}
We say that\/ $\Gamma:[0,T]\to D$ is {\em an accessible path\/}
 if\/ $\Gamma$ is a piecewise smooth curve composed from
a finite number of forward trajectories of the Hamiltonian systems
generated by $J_k$. 
\end{Definition}

More formally,
$\Gamma$ is an accessible path if there are $0=\tau_0<\tau_1<\dots<\tau_N=T$ such that
the sequence of points $z_i:=\Gamma(\tau_i)$ breaks the curve $\Gamma$ into trajectories, i.e.,
for every $i<N$, there is $k_i$, $1\le k_i \le n$, 
such that for $\tau\in [\tau_{i},\tau_{i+1}]$
\[
\Gamma(\tau)=\Phi_{k_i}^{\tau-\tau_{i} }(z_i)\,.
\]
Of course, the curve $\Gamma$ is well defined only if 
\[
0<\tau_{i+1}-\tau_{i}<\sigma_{k_i}(z_i)
\]
which ensures that the trajectories do not leave the domain $D$.

\begin{Theorem}\label{Thm:main}\label{TheoremSFH1}
If $D$ is a bounded domain in $\R^{2d}$,
the frozen fast system satisfies assumptions\/ {\rm [A1]} and {\rm [A2]}, 
$\{J_k\}_{k=1}^n$ is a set of actions corresponding
to a finite set of frozen periodic orbits in $\Lambda$,
and $\Gamma$ is an accessible path,
then there is a constant $C_0>0$ and $\varepsilon_0>0$
such that for every $\varepsilon<\varepsilon_0$
there is a trajectory of the full system  {\rm (\ref{Eq:motion})}
such that its slow component $z(t)$ satisfies
\[
 \|z(t)-\Gamma(\varepsilon t)\|<C_0 \varepsilon
\]
provided $0\le t \le \varepsilon^{-1}T$.
\end{Theorem}

\begin{Definition}
For any $z_0,z_1\in D$, we say that $z_1$ is {\em accessible\/}
from $z_0$ via the system $\{J_k\}$ if there is 
an accessible path such that $\Gamma(0)=z_0$ and $\Gamma(T)=z_1$.
\end{Definition}

In the case $d=1$ the accessibility property has a simple geometrical
meaning since trajectories of the Hamiltonian
systems generated by $J_k$ are level lines of the functions~$J_k$.
In this case the theorem provides trajectories which
follow segments of the level lines. The main obstacle for
the drift in the slow space is provided by level lines
common for all $J_k$.

\begin{Corollary}
Consider actions generated by two periodic orbits, $a$ and $b$. 
Those level lines of $J_{a,b}$, which are inside $D$,
are closed curves. The non-singular level lines form rings (or disks),
$D_a$ and $D_b$. Let $V=D_a\cap D_b\subset D$.
If $J_a$ and $J_b$ do not have common level lines, then
any point $z_1\in V$ is accessible from any point $z_0\in V$.
\end{Corollary}

\begin{Corollary}
Under the same assumptions.
Let us take any finite family of open sets $V_i\subset V$, which do not
depend on $\varepsilon$.
Then for all sufficiently small $\varepsilon$,
there is a trajectory which visits all the sets $V_i$. 
\end{Corollary}

If the energy set ${\cal M}_{\varepsilon}$ is compact
the slow dynamics never leaves a bounded set.
If at the same time $D$ is a connected set, 
natural questions arise: 
Is there a point in $D$ which 
is not accessible from every other point in $D$?
Is there a trajectory such that
its slow component is dense in $D$?

\section{Actions and first return maps near periodic orbits
of the frozen system\label{Se:action}}

Now consider the cylinder formed by periodic orbits of the frozen system:
\begin{equation}
S_{c,0}=\bigcup_{(v,u)\in D} L_c(v,u)\subset {\cal M}_0\,.
\end{equation}
 Let $\gamma_\varepsilon$ denote a trajectory of the full system (\ref{Eq:motion})
and $\pi_s:\R^{2m+2d}\to\R^{2d}$ the projection on the slow variables.

In the next section we will prove that some trajectories stay
in a neighbourhood of $S_{c,0}$ for a very long time and provide
a detailed description for them. In this section we show that 
in this case the evolution of the slow component  $\pi_s\gamma_\varepsilon$ 
approximately follows a trajectory of
the slow Hamiltonian flow $\Phi_c^{\varepsilon t}$ generated by the action $J_c$.

\begin{Lemma}\label{Le:stab}
Let $L_c$ be a family of periodic orbits of the frozen system.
If $\gamma_\varepsilon$ is a family of solutions of the full system\/~{\rm (\ref{Eq:motion})}
such that
\begin{itemize}
\item[{\rm (i)}]
 there are $z_0\in D$ and $C_0>0$ such that 
\begin{equation}\label{Eq:z0close}
\|\pi_s\gamma_\varepsilon(0)-z_0\|<C_0\varepsilon,
\end{equation}
\item[{\rm (ii)}]
 there are constants $C_1>0$ and $\tau_0<\sigma_c(z_0)$ such that
\begin{equation}\label{Eq:gammaclose}
\mathrm{dist}(\gamma_\varepsilon(t),S_{c,0})\le C_1\varepsilon\qquad
\forall t\in[0,\varepsilon^{-1}\tau_0],
\end{equation}
\end{itemize}
then there is $C_2>0$ such that
\begin{equation}\label{Eq:zdyn}
\| \pi_s\gamma_\varepsilon(t) - \Phi_c^{\varepsilon t}(z_0)\|\le C_2\varepsilon
\end{equation}
for all $t\in[0,\varepsilon^{-1}\tau_0]$.
\end{Lemma}

\noindent{\em Proof.}
We write
\(
(p_c(t,v,u),q_c(t,v,u))
\)
to denote a periodic solution of the frozen system
and use $T_c(v,u)$ for the corresponding period:
\begin{equation}
\label{Eq:per}
\begin{array}{rcl}
p_c(t+T_c(v,u),v,u)&\equiv& p_c(t,v,u)\,,\\
q_c(t+T_c(v,u),v,u)&\equiv &q_c(t,v,u)
\,.
\end{array}
\end{equation}
Then the action of the periodic orbit $L_c$
is given by the following integral
\begin{equation}\label{Eq:Jc}
J_c(v,u)
=\int_0^{T_c}p_c\frac{\partial q_c}{\partial t}\,dt\,.
\end{equation}
Since $L_c$ belongs to the zero energy level we have
a useful identity:
\begin{equation}\label{Eq:eplevel}
H(p_c(t,v,u),q_c(t,v,u),v,u;0)=0
\end{equation}
for all $(v,u)\in D$ and all $t\in\R$.

Let $\Sigma$ denote a smooth hypersurface  in $\R^{2m+2d}$ transversal to the flow 
of the frozen system such that  every periodic orbit of the family 
$L_c$ has exactly one intersection with~$\Sigma$. 
Let $M_i=\gamma_\varepsilon(t_i)$ be a sequence 
of consecutive intersections of $\gamma_\varepsilon$ with $\Sigma$
and consider the slow components of those points: $\hat z_i:=\pi_s M_i$.

We note that inequality (\ref{Eq:gammaclose}) 
and the smooth dependence of $p_c(s, z),q_c(s, z)$
on $z$ imply that there is $C_3>0$
such that for every $i$ there is $s_i$ such that
\[
\|M_i-(p_c(s_i,\hat z_i),q_c(s_i,\hat z_i),\hat z_i)\|\le C_1\varepsilon\,.
\]
Since solutions of differential equations
depend smoothly on the initial conditions and vector field,
the segment of $\gamma_\varepsilon(t)$, $t_i\le t\le t_{i+1}$ is 
close to $L_c(\hat z_i)$:
\begin{equation}\label{Eq:gap}
\gamma_\varepsilon(t)=(p_c(s_i+t-t_i,\hat z_i),q_c(s_i+t-t_i,\hat z_i),\hat z_i) + O(\varepsilon)
\end{equation}
and the time of the first return to the section $\Sigma$ is close to
the period of the frozen trajectory: 
\begin{equation}\label{Eq:tap}
t_{i+1}-t_i=T_c(\hat z_i)+O(\varepsilon).
\end{equation}
Now we estimate the displacement $\hat z_{i+1}-\hat z_i$.
We write $\hat z_i=(v,u)$ and $\hat z_{i+1}=(\bar v,\bar u)$ to shorten the notation. 
Integrating the slow component of the vector field along the exact trajectory
and using (\ref{Eq:motion}) we conclude
\begin{equation}\label{Eq:dudv}
\begin{array}{rclcl}
\bar u -u&=&\displaystyle
\int_{t_i}^{t_{i+1}} \dot u dt&=&\displaystyle
\varepsilon \int_0^{T_c(v,u)}
\left.\frac{\partial H}{\partial v}\right|_{p_c(t,v,u),q_c(t,v,u),v,u} dt+O(\varepsilon^2)\,,
\\
\bar v -v&=&\displaystyle
\int_{t_i}^{t_{i+1}} \dot v dt&=&\displaystyle
-\varepsilon \int_0^{T_c(v,u)}
\left.\frac{\partial H}{\partial u}\right|_{p_c(t,v,u),q_c(t,v,u),v,u} dt+O(\varepsilon^2)\,,
\end{array}
\end{equation}
where the error terms come from replacing the exact trajectory by the frozen one
and from the difference in the return time, see (\ref{Eq:gap}) and (\ref{Eq:tap}).
The integrals in the right hand side can be expressed in terms of derivatives of the action
defined by integral (\ref{Eq:Jc}).
Indeed, differentiating (\ref{Eq:Jc}) with respect to $u$, integrating by parts
and taking into account (\ref{Eq:per}), we get
\begin{eqnarray*}
\frac{\partial J_c}{\partial u}&=&
\int_0^{T_c}\left(
\frac{\partial p_c}{\partial u}\frac{\partial q_c}{\partial t}
-
\frac{\partial q_c}{\partial u}\frac{\partial p_c}{\partial t}
\right)dt
\\
&=&
\int_0^{T_c}\left(
\frac{\partial p_c}{\partial u}\frac{\partial H}{\partial p}
+
\frac{\partial q_c}{\partial u}\frac{\partial H}{\partial q}
\right)dt\,.
\end{eqnarray*}
Then differentiating identity (\ref{Eq:eplevel}) we get
\[
\frac{\partial p_c}{\partial u}\frac{\partial H}{\partial p}
+
\frac{\partial q_c}{\partial u}\frac{\partial H}{\partial q}
=-\frac{\partial H}{\partial u}\,,
\]
where the derivatives are evaluated at $(p_c(t,v,u),q_c(t,v,u),v,u)$.
Consequently
\[
\frac{\partial J_c}{\partial u}=-
\int_0^{T_c(v,u)}
\left.\frac{\partial H}{\partial u}\right|_{p_c(t,v,u),q_c(t,v,u),v,u} dt\,.
\]
Repeating these arguments with $u$ replaced by $v$ we also get
\[
\frac{\partial J_c}{\partial v}=-
\int_0^{T_c(v,u)}
\left.\frac{\partial H}{\partial v}\right|_{p_c(t,v,u),q_c(t,v,u),v,u} dt\,.
\]
Substituting the last two equalities into (\ref{Eq:dudv}) we arrive to
\begin{equation}
\bar u =u-\varepsilon \frac{\partial J_c}{\partial v}+O(\varepsilon^2)\,,\qquad
\bar v =v+\varepsilon \frac{\partial J_c}{\partial u}+O(\varepsilon^2)\,.
\end{equation}
We see that the displacement between two consecutive intersections of $\gamma_\varepsilon$ with
section $\Sigma$ is approximated by the time-$\varepsilon T_c$ shift 
along a trajectory of the Hamiltonian vector field (\ref{Eq:slow})
generated by the Hamiltonian function $J_c$:
\[
\hat z_{i+1}=\hat z_i+\Phi_c^{\varepsilon T_c}(\hat z_i)+O(\varepsilon^2)\,.
\]
Inequality (\ref{Eq:z0close}) implies that $\hat z_0=z_0+O(\varepsilon)$.
Then a rather standard stability estimate can be used to show 
\begin{equation}
\hat z_i=\Phi_c^{i\varepsilon T_c}(z_0)+O(\varepsilon)
\qquad \mbox{for
$0\le i \le \varepsilon^{-1}\tau_0$.}
\end{equation}
Finally, we note that $\hat z_i=\pi_s\gamma_\varepsilon(t_i)$ 
where $t_i=i T_c+O(\varepsilon i)$ due to (\ref{Eq:tap}).
Between intersections with $\Sigma$ the slow component $\pi_s\gamma_\varepsilon$
changes by a value of the order of~$O(\varepsilon)$. Therefore
the estimate is extendable to values of $t$ between the intersections
and (\ref{Eq:zdyn}) follows immediately.
\qed

\medskip

We note that $J_c$ is preserved by $\Phi_c^\tau$ and therefore
$J_c\circ\pi_s$ is an adiabatic invariant for the restriction of the full dynamics 
on a neighbourhood of $S_{c,0}$.

In general, we do not expect the estimates to be valid
on time intervals longer than stated  by Lemma~\ref{Le:stab}. 
For example, note that a trajectory 
of $J_c$ may leave the domain $D$ in finite time.

It is interesting that under additional assumptions
$\hat J_c(t):=J_c(\pi_s\gamma_\varepsilon(t))$ may stay
near its initial value, $\hat J_c(0)$, for much longer time.

Indeed, consider the case of one slow degree of freedom, $d=1$,
and suppose that level lines of $J_c$ are closed curves
on the plane of $(v,u)$ variables. Then $\sigma_c(z)=+\infty$
for all $z\in D$. In the next section we will show that the full system
has a locally invariant cylinder $S_{c,\varepsilon}$ close to $S_{c,0}$. 
Then equations (\ref{Eq:dudv}) describe a Poincar\'e map on the section 
defined by intersection of $S_{c,\varepsilon}$ and $\Sigma$.
In the case of one slow degree of freedom we may suppose that  the
map  (\ref{Eq:dudv}) satisfies assumptions of the KAM theorem,
then $\hat J_c(t)$ will stay close to its initial value forever.
Indeed, under these assumptions the trajectories on $S_{c,\varepsilon}$ 
are  trapped between two KAM tori. 

We also note that averaging theory can be used to study the
dynamics on $S_{c,\varepsilon}$ for $d\ge1$.

\section{Dynamics of the frozen system and normal hyperbolicity
\label{Se:normalhyperb}}

The following arguments are based on \cite{GT2007}.
Let $w=(p,q)$ and $z=(v,u)$ to shorten notation. 
Then the frozen system (\ref{Eq:frozen}) has the form
\begin{equation}\label{aham}
\dot w=G(w,z)\,,
\end{equation}
where $G$ is expressed in terms of partial derivatives of $H$
for $\varepsilon=0$.
The Hamiltonian function $H(w,z)$ is an integral of system (\ref{aham}).

Let system (\ref{aham}) have two smooth families of saddle periodic orbits
$L_a: w=w_a(t,z)$ and $L_b: w=w_b(t,z)$ for all $z\in D$.
Assume that both families belong to the zero energy level ${\cal M}_0$.
Take a pair of smooth cross-sections, $\Sigma_a$ and $\Sigma_b$,
which are transverse to the vector field and such that
each periodic trajectory $L_a(z)$ and $L_b(z)$ has exactly
one point of intersection with the corresponding section.
Denote the Poincar\'e map on $\Sigma_c$ near $L_c$ as $\Pi_{cc}$ ($c=a,b$).
The Poincar\'e map is smooth and depends smoothly on $z$.

We assume that for all $z\in D$ the frozen system has a pair of 
transversal heteroclinic orbits:
$\Gamma_{ab}\subseteq W^u(L_a)\cap W^s(L_b)$ and $\Gamma_{ba}\subseteq W^u(L_b)\cap W^s(L_a)$.

Let $\Pi_{ab}$ and $\Pi_{ba}$ be maps defined on subsets of $\Sigma_a$ and $\Sigma_b$
by following orbits close to $\Gamma_{ab}$ and $\Gamma_{ba}$, respectively.
Therefore $\Pi_{ab}$ acts from some open set in $\Sigma_a$ into an open set in $\Sigma_b$,
 while $\Pi_{ba}$ acts from an open set in $\Sigma_b$ into an open set in $\Sigma_a$. 
There is a certain freedom in the definition of the maps $\Pi_{ab}$ and $\Pi_{ba}$.
Indeed, each of these maps acts from a neighbourhood of one point of a heteroclinic 
orbit to a neighbourhood of another point of the same orbit, therefore different 
choices of the points lead to different maps.

\begin{figure}[ht]
\begin{center}
\includegraphics[width=0.5 \textwidth]{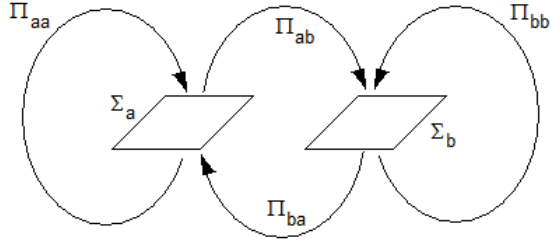}
\end{center}

\caption{Poincar\'e maps near two periodic orbits}
\end{figure}

When the maps are fixed, every orbit that lies entirely in a sufficiently small
neighbourhood of the heteroclinic cycle $L_a\cup L_b\cup \Gamma_{ab}\cup\Gamma_{ba}$
corresponds to a uniquely defined sequence of points $M_i\in\Sigma_a\cup\Sigma_b$ 
such that
$$
M_{i+1}=\Pi_{\xi_i\xi_{i+1}} M_i
$$
where
$$
\xi_i=c \;\mbox{ if }\; M_i\in\Sigma_c \;\;(c\in\{\,a,b\,\}).
$$
In this way the trajectory of the frozen system defines a 
sequence $\{\xi_i\}_{i=-\infty}^{i=+\infty}$
which is called {\em the code of the orbit}.

The periodic orbits $L_a$ and $L_b$ are saddle and the intersections
of the stable and unstable manifolds of $L_a$ and $L_b$ that create 
the heteroclinic orbits are transverse due to the assumptions [A1] and~[A2].
Consequently (cf. \cite{AfShi73}), one can choose the maps $\Pi_{ab}$ and $\Pi_{ba}$ 
and define coordinates $(x,y,z)$ in $\Sigma_a$ and  $\Sigma_b$ in such a way 
that the following holds.
\begin{itemize}
\item[{[H1]}]
$\Sigma_c\cap{\cal M}_0$ is diffeomorphic to the product $X_c\times Y_c\times D$ where $X_{c}$ and $Y_{c}$
are balls in $\mathbb{R}^{m-1}$ of a radius $R>0$ centred around the origin.

\item[{[H2]}] For each pair $c,c'\in\{\,a,b\,\}$ the Poincar\'e map $\Pi_{cc'}$ can be written in
the following ``cross-form'' \cite{book}: there exist smooth functions
$f_{cc'},g_{cc'}: X_c\times Y_{c'}\rightarrow X_{c'}\times Y_{c}$
such that a point $M(x,y,z)\in \Sigma_c$ is mapped
to $\bar M(\bar u,\bar w,z)\in \Sigma_{c'}$ 
by the map $\Pi_{cc'}$ if and only if
\begin{equation}\label{crossmap}
\bar x = f_{cc'}(x,\bar y,z), \qquad y = g_{cc'}(x,\bar y,z).
\end{equation}

\item[{[H3]}] There exists $\lambda<1$ such that
\begin{equation}\label{hyperb}
\left\|\frac{\partial (f_{cc'},g_{cc'})}{\partial (x,\bar y)}
\right\| \leq \lambda <1\,,
\end{equation}
where the norm of the Jacobian matrix corresponds to $\max\{\|x\|,\|y\|\}$.
\end{itemize}

Inequality (\ref{hyperb}) implies that for a fixed $z\in D$
the set $\Lambda_{z}$ of all orbits that
lie entirely in a sufficiently small neighbourhood of the heteroclinic
cycle $L_a\cup L_b\cup \Gamma_{ab}\cup\Gamma_{ba}$ in the energy level ${\cal M}_0$
 is hyperbolic, a horseshoe.
Moreover, one can show that the orbits in $\Lambda_{z}$ are in one-to-one correspondence with the set
of all sequences of $a$'s and $b$'s, i.e. for every sequence
$\{\xi_i\}_{i=-\infty}^{i=+\infty}$ there exists one and only one orbit
in $\Lambda_{z}$ which has this sequence as its code. 

Indeed, take any orbit from $\Lambda_{z}$ and denote by 
$M_i(x_i,y_i,z)$ the sequence of its intersections with 
the cross sections. Equation (\ref{crossmap}), 
implies that  the orbit has a code $\{\xi_i\}_{i=-\infty}^{i=+\infty}$
if and only if the coordinates of $M_i$ satisfy the equations
$$
x_{i+1}=f_{\xi_i\xi_{i+1}}(x_i,y_{i+1},z),
\qquad y_i = g_{\xi_i\xi_{i+1}}(x_i,y_{i+1},z)\,.
$$
Therefore the sequence $\{\,(x_i,y_i)\,\}_{i=-\infty}^{+\infty}$ is a fixed point of the operator
\begin{equation}\label{Eq:controp}
\{\,(x_i,y_i)\,\}_{i=-\infty}^{+\infty}\mapsto \{\,(f_{\xi_{i-1}\xi_i}(x_{i-1},y_i,z),
g_{\xi_i\xi_{i+1}}(x_i,y_{i+1},z)\,\}_{i=-\infty}^{+\infty}\,.
\end{equation}
Equation (\ref{hyperb}) implies this operator is a contraction of the space
$\prod_{i=-\infty}^{+\infty} X_{\xi_i}\times Y_{\xi_i}$, hence the
existence and uniqueness of the orbit with the code $\{\xi_i\}_{i=-\infty}^{i=+\infty}$
follow (see e.g. \cite{Shi67}). 

Moreover, the fixed point of a smooth contracting 
map depends smoothly on parameters. Consequently
the orbit depends smoothly on $z$ and the derivatives of
$(x_i(z,\xi),y_i(z,\xi))$ are bounded uniformly for all $i$.

\begin{Lemma}
If the Poincar\'e maps satisfy assumptions\/ {\rm [H1]--[H3]},
then for any two code sequences 
$\xi^{(1)}=\{\xi^{(1)}_i\}_{i=-\infty}^{+\infty}$ 
and $\xi^{(2)}=\{\xi^{(2)}_i\}_{i=-\infty}^{+\infty}$,
which satisfy
\[
\xi^{(1)}_i=\xi^{(2)}_i
\qquad \mbox{for $|i|\leq n$}
\]
the corresponding intersections with the cross sections
are bounded by
\begin{equation}\label{mix}
\max\left\{\,
\bigl\|x_{i}(z,\xi^{(1)})-x_{i}(z,\xi^{(2)})\bigr\|,\
\bigl\| y_i(z,\xi^{(1)})-y_i(z,\xi^{(2)})\bigr\|\,\right\}
\leq 2R \lambda^{n-|i|}\,,
\end{equation}
where the constants $R>0$ and $\lambda\in(0,1)$ are defined in\/ {\rm [H1]} and\/~{\rm [H3]}
respectively and do not depend on the sequences $\xi^{(1,2)}$.
\end{Lemma}

\noindent{\em Proof.} First we note, that 
\begin{eqnarray}\label{Eq:contr}
\lefteqn{ \max\left\{\,
\bigl\|x_{i}(z,\xi^{(1)})-x_{i}(z,\xi^{(2)})\bigr\|,\
\bigl\| y_{i+1}(z,\xi^{(1)})-y_{i+1}(z,\xi^{(2)})\bigr\|\,\right\}
}
\nonumber
\\
&\le&
\lambda
 \max\left\{\,
\bigl\|x_{i+1}(z,\xi^{(1)})-x_{i+1}(z,\xi^{(2)})\bigr\|,\
\bigl\| y_{i}(z,\xi^{(1)})-y_{i}(z,\xi^{(2)})\bigr\|\,\right\}
\end{eqnarray}
for $-n\le i<n$. 
Since none of the normes involved exceeds $2R$ 
we immediately conclude that
\begin{equation}\label{Eq:lr}
\max\left\{\,
\bigl\|x_{i}(z,\xi^{(1)})-x_{i}(z,\xi^{(2)})\bigr\|,\
\bigl\| y_{i+1}(z,\xi^{(1)})-y_{i+1}(z,\xi^{(2)})\bigr\|\,\right\}
\le
2 R\lambda\,.
\end{equation}
Then the following estimate is true for $n'=1$
\begin{eqnarray}\label{Eq:ind}
&\max\bigl\{\,
\bigl\|x_{i}(z,\xi^{(1)})-x_{i}(z,\xi^{(2)})\bigr\|,\
\bigl\| y_{i+1}(z,\xi^{(1)})-y_{i+1}(z,\xi^{(2)})\bigr\|\,\bigr\}
\nonumber
\\
&\quad\le
2 R\left\{
\begin{array}{ll}
\lambda^{n'-i}&\mbox{for $0\le i< n'$}\\
\lambda^{n'+1+i}&\mbox{for $-n'-1< i< 0$}
\end{array}
\right.
\,.
\end{eqnarray}
We continue inductively in $n'$.
Assuming that the estimate (\ref{Eq:ind}) holds for $n'$ replaced by $n'-1$
we check the upper bounds for all $2n'+1$ different values of $i$
in the order of decreasing of $|i|$. On each step we use the contraction
property (\ref{Eq:contr}) and the sharper of upper bounds (\ref{Eq:ind}) and (\ref{Eq:lr}).
In the case $n'=n$,   
\[
\max\bigl\{\,
\bigl\|x_{n}(z,\xi^{(1)})-x_{n}(z,\xi^{(2)})\bigr\|,\
\bigl\| y_{-n}(z,\xi^{(1)})-y_{-n}(z,\xi^{(2)})\bigr\|\,\bigr\}
\le
2 R
\]
is used instead of (\ref{Eq:lr}). Then (\ref{mix}) follows directly from 
the last upper bound and (\ref{Eq:ind}) taken with $n'=n$.
\qed

\medskip

Now let us consider the full system (\ref{Eq:motion}) for a small $\varepsilon>0$.
Since the vector field depends smoothly on $\varepsilon$,
the Poincar\'e maps $\Pi_{cc'}: \Sigma_c\rightarrow \Sigma_{c'}$
are still defined and can be written in the following form:
\begin{equation}\label{crosseps}
\left\{\begin{array}{ll}\displaystyle
\bar x = f_{cc'}(x,\bar y,z,\varepsilon), \qquad y = g_{cc'}(x,\bar y,z,\varepsilon)\\[6pt]
\displaystyle \bar z= z+\varepsilon \phi_{cc'}(x,\bar y,z,\varepsilon),
\end{array}\right.
\end{equation}
where $f,g,\phi$ are bounded along with the
first derivatives and $f,g$ satisfy (\ref{hyperb}).

For technical reasons we need to assume that the domain $D$ is invariant
under the Poincar\'e map, i.e., $\phi_{cc'}(x,\bar y,z,\varepsilon)\equiv 0$ 
if $z\in\partial D$.
If this is not the case, then we modify $\phi_{cc'}$ in a small neighbourhood of 
the boundary. 
We note that the next lemma contains a statement of uniqueness but
the surfaces provided by the lemma may depend on the way the functions $\phi_{cc'}$
have been modified. Therefore the lemma implies existence but not uniqueness
for the original system.

The next lemma is of a general nature and has little to do with the Hamiltonian
structure of the equations. Rather we notice that by fixing any code $\xi$ and varying $z\in D$
we obtain at $\varepsilon=0$ a sequence of smooth two-dimensional surfaces.
The $i$-th surface is the set run by the point $M_i(z)$ of the uniquely
defined orbit with the code $\xi$. This sequence is invariant with respect to the
corresponding Poincar\'e maps and is uniformly normally-hyperbolic --- hence it
persists at all $\varepsilon$ sufficiently small.

\begin{Lemma}\label{Lemma1}
Given any sequence $\xi$ of $a$'s and $b$'s,
there exists a uniquely defined sequence of smooth surfaces
\begin{equation}\label{seq}
{\cal L}_i(\xi,\varepsilon):
(x,y)=(x_i(z,\xi,\varepsilon),y_i(z,\xi,\varepsilon))
\end{equation}
such that
\begin{equation}\label{inva}
\Pi_{\xi_i\xi_{i+1}} {\cal L}_i={\cal L}_{i+1}.
\end{equation}
The functions $(x_i,y_i)$ are defined for all small $\varepsilon$ and all
$z\in D$, they are
uniformly bounded along with their derivatives with respect to $z$
and satisfy~{\rm (\ref{mix})}.
Moreover, there is $C>0$ independent from the code $\xi$ such that
\[
\bigl\|
(x_i(z,\xi,\varepsilon)-x_i(z,\xi,0),
y_i(z,\xi,\varepsilon)-y_i(z,\xi,0))\bigr\|
\le C\varepsilon\,,
\]
for all $i\in\mathbb{Z}$.
\end{Lemma}

A proof of this lemma is essentially identical to the proof of Lemma~1
of~\cite{GT2007} and is based on contraction mapping arguments:
the functions $x_i,y_i$ are constructed as a fixed point of
an operator similar to (\ref{Eq:controp}).

We note that if $\xi=c^\infty$ is a code which consists of the symbol $c$
only, then ${\cal L}_i$, $x_i$ and $y_i$ are independent from $i$,
and we will denote them by ${\cal L}_c$, $x_c$ and $y_c$ respectively.

\section{Drift of slow variables}

Let $\xi $ be a code. The corresponding trajectory of the full system
is described by the dynamics of its slow component:
\begin{equation}\label{Eq:ximap}
z_{i+1}=z_{i}+\varepsilon \phi _{\xi_{i}\xi_{i+1}}
(x_{i}(z_{i},\xi ,\varepsilon ),y_{i}(z_{i},\xi ,\varepsilon
),z_{i},\varepsilon )
\,.
\end{equation}
If $\xi=c^\infty$ the functions $x_i$ and $y_i$ do not depend on $i$, 
so we denoted them by $x_c$ and $y_c$ respectively.
Then the equation can be written in the form
\begin{equation}\label{Eq:cmap}
\bar z_{i+1}=\bar z_{i}+\varepsilon \phi _{cc}
(x_{c}(\bar z_{i},\varepsilon ),y_{c}(\bar z_{i} ,\varepsilon
),\bar z_{i},\varepsilon )
\,,
\end{equation}
where the bars over $z_i$ and $z_{i+1}$ are used to distinguish trajectories
of (\ref{Eq:cmap}) and (\ref{Eq:ximap}).

The next lemma estimates the difference between these two
slow dynamics for all sequences which have a large block of $c$'s.

\begin{lemma}
\label{LemmaSFH1_2}
Assume the assumptions of Lemma~\ref{Lemma1} are satisfied.
Then for any $K_0>0$, $t_0>0$, there is $K_1>0$ and $\varepsilon_0>0$ 
such that for any $|\varepsilon|<\varepsilon_0$
and any code $\xi$ such that for some index $j$
\[
\xi _{j}=\xi _{j+1}=\ldots =\xi _{j+\left\lfloor \frac{t_{0}}{\varepsilon }%
\right\rfloor }=c
\]
the inequality $\left\Vert z_{j}-\bar z_0 \right\Vert \leq \varepsilon K_{0}$ 
implies the corresponding trajectories of\/ {\rm (\ref{Eq:ximap})}
and\/  {\rm (\ref{Eq:cmap})} satisfy the inequality
\[
\left\Vert z_{j+N}^{{}}-\bar z_N\right\Vert \leq \varepsilon K_{1}
\] 
for all $0\leq N\leq N_0(\varepsilon)\equiv\left\lfloor \frac{t_{0}}{\varepsilon }\right\rfloor $.
\end{lemma}

\noindent{\em Proof.}
Using (\ref{Eq:ximap}) we get
\begin{equation}
z_{j+N}=z_{j}+\varepsilon \sum_{i=j}^{j+N-1}
\phi _{\xi _{i}\xi
_{i+1}}(x_{i}(z_{i},\xi ,\varepsilon ),y_{i}(z_{i},\xi ,\varepsilon
),z_{i},\varepsilon )\,. \label{SFH_Proof21}
\end{equation}%
Using (\ref{Eq:cmap}), 
we obtain in a similar way 
\begin{equation}
\bar{z}_{N}=\bar{z}_{0}+\varepsilon \sum_{i=0}^{N-1}
\phi_{cc}(x_{c}(\bar{z}_{i},\varepsilon),y_{c}(\bar{z}_{i},\varepsilon),\bar{z}_{i},\varepsilon).  
\label{SFH_Proof22}
\end{equation}%
We have assumed $\|\bar z_0-z_j\|\le K_0\varepsilon$. Then taking the difference of the equalities
(\ref{SFH_Proof21}) and (\ref{SFH_Proof22}), 
we obtain  
\begin{eqnarray}
\left\Vert z_{j+N}-\bar{z}_{N}\right\Vert  
&\leq & \varepsilon K_0+
\varepsilon \sum_{i=0}^{N-1}
\bigl\| 
\phi _{cc}(x_{j+i}(z_{j+i},\xi ,\varepsilon ),y_{j+i}(z_{j+i},\xi ,\varepsilon ),
z_{j+i},\varepsilon )   \nonumber \\
&&\qquad\qquad
-\phi _{cc}(x_{c}(\bar z_{i},\varepsilon ),y_{c}(\bar z_{i},\varepsilon ),\bar z_{i},\varepsilon )
\bigr\|
\label{SFH_Proof23}
\end{eqnarray}
Consequently,
\begin{eqnarray}
\label{SFH_Proof24}
\lefteqn{
\left\Vert z_{j+N}-\bar{z}_{N}\right\Vert \leq 
\varepsilon K_0+
\varepsilon \left\Vert \frac{\partial \phi _{cc}}{\partial z}\right\Vert 
\sum_{i=0}^{N-1}\left\| z_{j+i}-\bar{z}_{i}\right\|
}
 \\
&&+
\varepsilon \left\Vert \frac{\partial \phi _{cc}}{\partial (x,y)}%
\right\| 
\sum_{i=0}^{N-1}
\bigl\Vert 
\bigl( x_{j+i}(z_{j+i},\xi ,\varepsilon)-x_{c}(\bar{z}_{i},\varepsilon),
       y_{j+i}(z_{j+i},\xi ,\varepsilon )-y_{c}(\bar{z}_{i},\varepsilon)\bigr) 
\bigr\| .  
\nonumber
\end{eqnarray}%
In order to estimate the last term we note that Lemma~\ref{Lemma1} includes
the estimate~(\ref{mix})
\begin{eqnarray*}
\lefteqn{
\sum_{i=0}^{N-1}
\bigl\Vert 
\bigl( x_{j+i}(z_{j+i},\xi ,\varepsilon)-x_{c}(\bar{z}_{i},\varepsilon),
       y_{j+i}(z_{j+i},\xi ,\varepsilon )-y_{c}(\bar{z}_{i},\varepsilon)\bigr) 
\bigr\|
}
\\
&\le&
\sum_{i=0}^{N-1}
\bigl\Vert 
\bigl( x_{j+i}(z_{j+i},\xi ,\varepsilon)-x_{c}({z}_{j+i},\varepsilon),
       y_{j+i}(z_{j+i},\xi ,\varepsilon )-y_{c}({z}_{j+i},\varepsilon)\bigr) 
\bigr\|
\\
&&+\sum_{i=0}^{N-1}
\bigl\Vert 
\bigl( x_{c}(z_{j+i},\varepsilon)-x_{c}(\bar{z}_{i},\varepsilon),
       y_{c}(z_{j+i} ,\varepsilon )-y_{c}(\bar{z}_{i},\varepsilon)\bigr) 
\bigr\|
\\
&\le&
\sum_{i=0}^{N-1} 2R \lambda^{\min\{i,N_0(\varepsilon)-i\}}
+
\sum_{i=0}^{N-1}
\max\left\{\left\|
\frac{\partial x_{c}}{\partial z}
\right\|
,
\left\|
\frac{\partial y_{c}}{\partial z}
\right\|
\right\}
\|z_{j+i}-\bar{z}_{i}\|
\\
&\le&
\frac{4R}{1-\lambda}
+
\max\left\{\left\|
\frac{\partial x_{c}}{\partial z}
\right\|
,
\left\|
\frac{\partial y_{c}}{\partial z}
\right\|
\right\}
\sum_{i=0}^{N-1}
\|z_{j+i}-\bar{z}_{i}\|
\end{eqnarray*}%
Substituting the last bound into (\ref{SFH_Proof24})
we obtain
\begin{eqnarray*}
\left\Vert z_{j+N}-\bar{z}_{N}\right\Vert &\leq &
\varepsilon 
\left(
\left\Vert \frac{\partial \phi _{cc}}{\partial z}\right\Vert 
+
 \left\Vert \frac{\partial \phi _{cc}}{\partial (x,y)}\right\Vert
\max\left\{\left\|
\frac{\partial x_{c}}{\partial z}
\right\|
,
\left\|
\frac{\partial y_{c}}{\partial z}
\right\|
\right\}
\right)
\sum_{i=0}^{N-1}\left\| z_{j+i}-\bar{z}_{i}\right\|
\\
&&
+
\varepsilon \left\Vert \frac{\partial \phi _{cc}}{\partial (x,y)}\right\Vert
\frac{4R}{1-\lambda}+\varepsilon K_0\,.
\end{eqnarray*}
Let
\begin{eqnarray*}
A&=&
\left\Vert \frac{\partial \phi _{cc}}{\partial (x,y)}\right\Vert
\frac{4R}{1-\lambda}+ K_0\,.
\\
B&=&
\left\Vert \frac{\partial \phi _{cc}}{\partial z}\right\Vert 
+
\left\Vert \frac{\partial \phi _{cc}}{\partial (x,y)}\right\Vert
\max\left\{\left\|
\frac{\partial x_{c}}{\partial z}
\right\|
,
\left\|
\frac{\partial y_{c}}{\partial z}
\right\|\right\}\,.
\end{eqnarray*}
Then
\[
\left\Vert z_{j+N}-\bar{z}_{N}\right\Vert \leq 
\varepsilon A+\varepsilon B
\sum_{i=0}^{N-1}\left\| z_{j+i}-\bar{z}_{i}\right\|
\]
and consequently
\[
\left\Vert z_{j+N}-\bar{z}_{N}\right\Vert \leq 
\varepsilon A{\mathrm e}^{\varepsilon NB}\,.
\]
So for $N\leq \frac{t_{0}}{\varepsilon }$ we have%
\begin{equation}
\left\Vert z_{j+N}-\bar{z}_{N}\right\Vert \leq  K_{1}\varepsilon 
\label{SFH_Proof29}
\end{equation}
where $K_1=A{\mathrm e}^{t_0B}$.
\qed

\medskip

We note that Lemma~\ref{LemmaSFH1_2} is also valid for any two sequences
with any large common block.

\begin{lemma}\label{Le:any2}
Assume the assumptions of Lemma~\ref{Lemma1} are satisfied.
Then for any $K_0>0$, $t_0>0$, there is $K_1>0$ and $\varepsilon_0>0$ 
such that for any $|\varepsilon|<\varepsilon_0$
and any two codes $\xi^{(1)}$ and $\xi^{(2)}$ such that for some index $j$
\[
\xi^{(1)}_{j+i}=\xi^{(2)}_{j+i}
\qquad
0\le i\le N_0(\varepsilon)\equiv\left\lfloor \frac{t_{0}}{\varepsilon }\right\rfloor 
\]
the inequality $\left\Vert z_{j}^{(1)}- z_j^{(2)} \right\Vert \leq \varepsilon K_{0}$ 
implies
\[
\left\Vert z_{j+N}^{(1)}-z_{j+N}^{(2)}\right\Vert \leq \varepsilon K_{1}
\] 
for all $0\leq N\leq N_0(\varepsilon)$.
\end{lemma}

The proof of this lemma is almost identical to the previous one so we skip it.

\section{Proof of Theorem~\ref{Thm:main}}

Now we have all ingredients necessary 
to complete the proof of Theorem~\ref{Thm:main}.
Each periodic orbit $L_k\in\Lambda$ is defined by 
a periodic code. Let $\ell_0$ be the longest period
among the codes corresponding to the periodic orbits
selected in the assumptions of Theorem~\ref{Thm:main}.
Then each of the periodic orbits can be uniquely
identified by a piece of code $c_k\in\{\,a,b\,\}^{\ell_0}$.

Given an accessible path $\Gamma$ we define
\[
\Delta_i=\tau_{i+1}-\tau_{i}\,.
\]
It is the time the slow motion follows the flow defined
by the Hamiltonian function $J_{k_i}$, $1\le i\le N$,
where $N$ is the number of segments in the path. 
Let 
\[
N_i(\varepsilon)=\left\lfloor \frac{\Delta_i}{\varepsilon \ell_0}\right\rfloor\,.
\]
Let $\omega_i=c_{k_i}^{N_i(\varepsilon)}$ be a finite sequence, which
consists of $N_i(\varepsilon)$ copies of the symbol $c_{k_i}$.
Let $\xi_\varepsilon$ be any sequence, which contains $\omega_1\omega_2\dots\omega_N$
starting from position $0$. Obviously, the sequence 
\[
j_i=\ell_0\sum_{l=0}^i N_l(\varepsilon)
\]
indicates the starting positions of the blocks $\omega_i$ in the code $\xi_\varepsilon$.

We note that assumptions [A1] and [A2] imply that there are sections
and Poincar\'e maps of the frozen system (\ref{Eq:frozen}) 
which satisfy [H1]---[H3].
In order to apply Lemma~\ref{Lemma1} we have to modify the slow
component of the Poincar\'e maps to ensure invariance of $D$.
Since $D$ is open there is $\delta>0$ such that a $\delta$-neighbourhood
of $\Gamma$ is inside $D$. Then we modify $\phi_{cc'}$ outside 
this $\delta$-neighbourhood of $\Gamma$ to ensure that $\phi_{cc'}$
vanishes near $\partial D$. This modification does
not affect trajectories which do not leave
an $O(\varepsilon)$ neighbourhood of $\Gamma$: i.e. while a trajectory of
the modified maps stays in the neighbourhood of $\Gamma$ it is simultaneously
a trajectory of the original Poincar\'e maps.

Lemma~\ref{Lemma1} implies that there is a sequence of surfaces which corresponds to
the sequence $\xi\equiv\xi_\varepsilon$. Now let $z_0=\Gamma(0)$ and consider
the sequence of points 
\[
M_i\equiv
M_{\xi_i}(x_i(z_i,\xi,\varepsilon),y_i(z_i,\xi,\varepsilon),z_i;\varepsilon)
\]
on those surfaces. The slow component $z_i$ satisfies equation (\ref{Eq:ximap})
and
\[
M_i\in {\cal M}_\varepsilon\cap \Sigma_{\xi_i }\,.
\]
Lemma~\ref{LemmaSFH1_2} (or Lemma~\ref{Le:any2}) implies that
\[
\|z_i-\bar z_i\|\le K_1 \varepsilon
\]
for all $i$ such that $0\le i \le N_1(\varepsilon)$, where $\bar z_i$ 
denote the trajectory of (\ref{Eq:cmap}). 
We continue inductively to show using Lemma~\ref{LemmaSFH1_2} 
that there are constant $K_k$ such that
\begin{equation}
\|z_i-\bar z_{i-j_{k-1}}\|\le K_k \varepsilon
\end{equation}
for all $i$ such that $j_{k-1}\le i\le j_k$ and $1\le k\le N$,
where $\bar z_l$ satisfy (\ref{Eq:ximap}) with initial condition $\bar z_0:=z_{j_{k-1}}$.
We note that these $\bar z_{i-j_{k-1}}$ all belong to the invariant surface ${\cal L}_{c_k}$,
then Lemma~\ref{Lemma1} implies 
\[
\mathrm{dist}(M_i,S_{c_k,0})=O(\varepsilon)\,.
\]
We consider the trajectory of the full system (\ref{Eq:motion}),
which we denote by $\gamma_\varepsilon$, such that $\gamma_\varepsilon(0)=M_0$.
Since $\gamma_\varepsilon$ goes through the points $M_i$
it also stays $O(\varepsilon)$-close to $S_{c_k,0}$ between the points
therefore
\[
\mathrm{dist}(\gamma_\varepsilon(t),S_{c_k,0})=O(\varepsilon)
\]
for $t\in[\tau_{k-1}\varepsilon^{-1},\tau_k\varepsilon^{-1}]$.
Then Lemma~{\ref{Le:stab}} implies that
the slow component $\pi_s\gamma_\varepsilon(t)$ 
shadows the accessible path $\Gamma$.

\section{Final remarks}
Finally, we note that similar equations arise in the case of a Hamiltonian system
 with the standard symplectic form, $\Omega_{\mathrm{st}}=dp\wedge dq +dv \wedge du$ 
when a Hamiltonian function looses some degrees of freedom at $\varepsilon=0$. 
More precisely, if the Hamiltonian function has the form
$$
{\cal H}(p,q,v,u;\varepsilon)=H_0(p,q)+\varepsilon H_1(p,q,v,u)\,,
$$
the corresponding Hamilton equations are given by
\begin{equation}\label{Eq:motion1}
\begin{array}{ll}
\displaystyle
\dot q=\frac{\partial H_0}{\partial p}+\varepsilon \frac{\partial H_1}{\partial p}\,, \qquad&
\displaystyle
\dot p=-\frac{\partial H_0}{\partial q}-\varepsilon \frac{\partial H_1}{\partial q}\,,\\[10pt]
\displaystyle
\dot u=\varepsilon\frac{\partial H_1}{\partial v}\,, &
\displaystyle
\dot v=-\varepsilon\frac{\partial H_1}{\partial u}  \,.
\end{array}
\end{equation}
In this equation, adiabatic invariants can be
destructed by resonances~\cite{NeishtadtVasiliev}.

These equation are quite similar to (\ref{Eq:motion}).
We note that the frozen fast system is independent of
the slow variables. The theory developed in this
paper can be applied to the system (\ref{Eq:motion1}).
The most notable difference is related to the description
of the slow motion near a cylinder formed by periodic 
orbits of the frozen system. Indeed the slow motion
is described  by the averaged perturbation term
\[
\tilde J_c(u,v)=
\int_0^{T_c} H_1(p_c(t),q_c(t),u,v)\,dt
\]
and not by the actions.


\begin{thebibliography}{GL2007}
\bibliographystyle{alpha}

\bibitem{Arnold78}
Arnold V.I. Mathematical Methods of Classical Mechanics. Springer Verlag, 1978

\bibitem{BM1961}
Bogolyubov N.N., Mitropol'skii Yu.A. Asymptotic
Methods in the Theory of Nonlinear Oscillations
Gordon and Breach, 1961.


\bibitem{AfShi73} Afraimovich, V.S., Shilnikov, L.P., On critical sets of Morse-Smale
systems, Trans. Moscow Math. Soc. 28 (1973) 179--212.

\bibitem{Anosov}
Anosov D., Averaging in systems of ODEs with rapidly oscillating solutions, 
Izv. Akad. Nauk. SSSR 24 (1960) 721--742

%
%

\bibitem{GelfreichLerman}
 Gelfreich V., Lerman L. 
    Long-periodic orbits and invariant tori in a singularly perturbed Hamiltonian system, 
    Physica D Vol 176 Iss. 3--4, (2003) 125--146

\bibitem{GT2007}
Gelfreich V., Turaev D.,
Unbounded energy growth in Hamiltonian systems with a slowly varying parameter,
Math. Physics Preprint Archive ({\tt http://www.ma.utexas.edu/mp\_arc}), 
preprint 07-215 (2007) 30 p. 

\bibitem{Kiefer}
Y. Kifer. Another proof of the averaging principle for fully coupled dynamical systems with hyperbolic fast
motions, Disc. and Cont. Dynam. Sys., Vol. 13, No. 5 (2005) 1187--1201.

\bibitem{Lochak}
P. Lochak and C. Meunier. Multiphase averaging for classical systems, Springer Verlag, New York, 1988.

\bibitem{Neishtadt1976}
Neishtadt A., Averaging in multi-frequency systems. II. Sov. Phys. Dokl, 21 (1976) 80--82.

\bibitem{NeishtadtVasiliev}
Neishtadt A., Vasiliev A.,
Destruction of adiabatic invariance at resonances in slow-fast Hamiltonian systems
Physics Research A 561, 158-165 (2006)

\bibitem{Shi67} Shilnikov, L.P., On a Poincar\'e-Birkhoff problem, Math. USSR Sb. 3 (1967) 91--102.

\bibitem{book} Shilnikov, L.P., Shilnikov, A.L., Turaev, D.V., Chua, L.O.,
Methods of qualitative theory in nonlinear dynamics. Part I. World Scientific, 1998.

\end{thebibliography}
\end{document}